\documentclass{amsart}

\usepackage{tikz}
\usepackage[margin=1in]{geometry}
\usepackage{amsmath}
\usepackage{amssymb}
\usepackage{amsthm}
\usepackage{cite}
\usepackage{hyperref}
\usepackage{algorithmic}
\usepackage{setspace}
\usepackage{mathtools}

\newtheorem{theorem}{Theorem}[section]
\newtheorem{lemma}[theorem]{Lemma}

\newtheorem{proposition}[theorem]{Proposition}

\author{Hunter Waldron}
\date{}
\title{A Combinatorial Proof of a Partition Perimeter Inequality}

\begin{document}

\begin{abstract}
The partition perimeter is a statistic defined to be one less than the sum of the
number of parts and the largest part. Recently,
Amdeberhan, Andrews, and Ballantine proved the following analog of Glaisher's theorem:
for all $m \geq 2$ and $n \geq 1$, there are at least as many partitions with
perimeter $n$ and parts $\not \equiv 0 \pmod{m}$ as partitions with perimeter $n$
and parts repeating fewer than $m$ times. In this work, we provide a combinatorial
proof of their theorem by relating the combinatorics of the partition perimeter to
that of compositions. Using this technique, we also show that a composition theorem
of Huang implies a refinement of another perimeter theorem of Fu and Tang.
\end{abstract}

\maketitle

\section{Introduction}

The \textit{partitions} of an integer $n \geq 1$, which are finite sequences of integers
$\lambda = (\lambda_1, \lambda_2, \dots, \lambda_\ell)$ where the \textit{parts}
$\lambda_i$ sum to $n$ and satisfy the inequalities
$\lambda_1 \geq \lambda_2 \geq \cdots \geq \lambda_\ell$, have long been studied
for their wealth of enumerative properties. Many seemingly unrelated families of
partitions have been shown to be equinumerous. Euler's theorem is a classic example
of this kind of result:
\begin{theorem}
For all $n \geq 1$, there are as many partitions of $n$ with odd parts as there are
partitions of $n$ with distinct parts.
\end{theorem}

The \textit{perimeter} of a partition is a statistic defined to be the value
$\lambda_1 + \ell(\lambda) - 1$
where $\ell(\lambda)$ is the number of parts in $\lambda$. This is
Straub's definition \cite{Straub}, which is the most convenient for the purposes of this
work. Many other partition statistics like the \textit{rank} $\lambda_1 - \ell(\lambda)$
allow infinitely many partitions to take on each value, and so must be studied in
addition to partition \textit{size} $n$. There are only $2^{k-1}$ partitions with perimeter
$k \geq 1$, which makes viable the study of partitions by their perimeter alone.

Interest in enumeration by perimeter
(see \cite{AmdeberhanAndrewsBallantine, FuTang, LinXiongYan})
spawned from the following result of Straub, which is
remarkably similar in presentation to Euler's theorem.

\begin{theorem}[Theorem 1.4 in \cite{Straub}] % TODO fix this
For all $n \geq 1$, there are as many partitions with perimeter $n$ and odd parts
as there are partitions with perimeter $n$ and distinct parts.
\end{theorem}

Naturally, for any theorem like Euler's we can attempt to replace size with perimeter
and show that an analogous result holds. For instance,
Glaisher's theorem generalizes Euler's to allow an arbitrary modulus:

\begin{theorem}
For all $m \geq 2$ and $n \geq 1$, there are as many partitions of $n$ with parts
$\not \equiv 0 \pmod{m}$ as there are partitions of $n$ with parts repeating fewer
than $m$ times.
\end{theorem}

Let $h_m(n)$ and $g_m(n)$ be the sets of partitions with perimeter $n$ having
parts $\not \equiv 0 \pmod{m}$ and parts repeating fewer than $m$ times,
respectively. In \cite{AmdeberhanAndrewsBallantine},
Amdeberhan, Andrews, and Ballantine showed that
Straub's theorem actually becomes an inequality when generalized in the direction
of Glaisher's theorem. 

\begin{theorem}[Theorem 7.6 in \cite{AmdeberhanAndrewsBallantine}]
\label{main_result}
We have $h_m(n) - g_m(n) \geq 0$ for all $m \geq 2$ and $n \geq 1$.
\end{theorem}
They proved this result by manipulating generating functions,
and asked if a combinatorial proof exists. In this work, we provide such a proof.
The technique we use appears to have additional utility in this kind of research.
As an example, we provide another application to prove a refinement of the following
theorem of Fu and Tang.

\begin{theorem}[Theorem 2.15 in \cite{FuTang}]
\label{theorem_fu_tang}
For all $m, n \geq 1$, there are as many partitions with parts $\equiv 1 \pmod{m+1}$
and perimeter $n$ as there are partitions with gaps between parts $\geq m$
and perimeter $n$.
\end{theorem}

Our refinement, which we state here,
adds the additional parameter $k$, and reduces to Fu and Tang's theorem
when $k=0$, and then to Straub's theorem when $m=1$:

\begin{theorem}
\label{theorem_fu_tang_generalization}
For all $n, m \geq 1$ and $k \geq 0$, the following families of partitions with
perimeter $n$ are equinumerous:
\begin{itemize}
\item[1)] Partitions $\lambda$ with parts
grouped into $k$ non-empty and disjoint sets of sequential parts $\{\lambda_1, \dots, \lambda_{t_1}\}$,
$\{\lambda_{t_1+1},$ $\dots,$ $\lambda_{t_2}\},$ $\dots,$
$\{\lambda_{t_{k-1}+1}, \dots,\lambda_{t_k = \ell(\lambda)}\}$, or $k+1$ sets if
$\lambda_{\ell(\lambda)} \equiv 1 \pmod{m+1}$ but with  $\lambda_{\ell(\lambda)} >
m+1$ otherwise, where each set has parts belonging all to the same congruence class
modulo $m+1$, and neighboring sets have different congruence classes and gaps
between them $\lambda_{t_i} - \lambda_{t_{i+1}} > m$ for each $1 \leq i < k$.

\item[2)] Partitions $\lambda$ with gaps between parts $\geq m$ with
exactly $k$ exceptions $\lambda_i - \lambda_{i+1} < m$,
always proceeded by a gap $\lambda_{i-1} - \lambda_i \geq m$, and
followed by a gap $\lambda_{i+1} - \lambda_{i+2} > m$ unless $\ell(\lambda)
= i+1$.
\end{itemize}
\end{theorem}

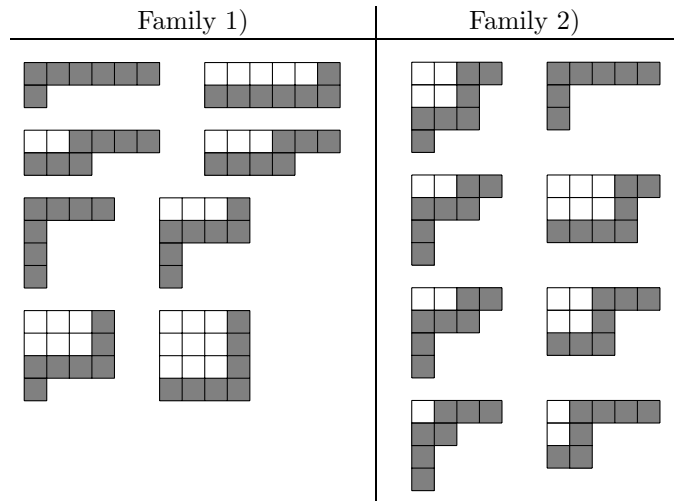
\begin{figure}[!ht]
\begin{tabular}{c | c}
Family 1) & Family 2) \\
\hline
\begin{tikzpicture}[scale=0.3]

\draw[draw=none] (14,0) -- (15,0);

\fill[color=gray] (0,0) rectangle (6,-1);
\fill[color=gray] (0,-1) rectangle (1,-2);
\draw (0,0) grid (6,-1);
\draw (0,-1) grid (1,-2);

\fill[color=gray] (13,0) rectangle (14,-1);
\fill[color=gray] (8,-1) rectangle (14,-2);
\draw (8,0) grid (14,-2);

\fill[color=gray] (2,-3) rectangle (6,-4);
\fill[color=gray] (0,-4) rectangle (3,-5);
\draw (0,-3) grid (6,-4);
\draw (0,-4) grid (3,-5);

\fill[color=gray] (11,-3) rectangle (14,-4);
\fill[color=gray] (8,-4) rectangle (12,-5);
\draw (8,-3) grid (14,-4);
\draw (8,-4) grid (12,-5);

\fill[color=gray] (0,-6) rectangle (4,-7);
\fill[color=gray] (0,-7) rectangle (1,-10);
\draw (0, -6) grid (4,-7);
\draw (0,-7) grid (1, -10);

\fill[color=gray] (9,-6) rectangle (10,-7);
\fill[color=gray] (6,-7) rectangle (10,-8);
\fill[color=gray] (6,-8) rectangle (7,-10);
\draw (6,-6) grid (10,-8);
\draw (6,-8) grid (7,-10);

\fill[color=gray] (3,-11) rectangle (4,-14);
\fill[color=gray] (0,-13) rectangle (3,-14);
\fill[color=gray] (0,-14) rectangle (1,-15);
\draw (0, -11) grid (4, -14);
\draw (0,-14) grid (1,-15);

\fill[color=gray] (9,-11) rectangle (10,-15);
\fill[color=gray] (6,-14) rectangle (9,-15);
\draw (6,-11) grid (10,-15);

\draw[draw=none] (0,-15) -- (0,-19);

\end{tikzpicture} &
\begin{tikzpicture}[scale=0.3]
\draw[draw=none] (0,0) -- (-1,0);
\draw[draw=none] (0,0) -- (0,1);

\fill[color=gray] (2,0) rectangle (4,-1);
\fill[color=gray] (2,-1) rectangle (3,-3);
\fill[color=gray] (0,-2) rectangle (2,-3);
\fill[color=gray] (0,-3) rectangle (1,-4);
\draw (0,0) grid (4,-1);
\draw (0,-1) grid (3,-3);
\draw (0,-3) grid (1,-4);

\fill[color=gray] (6,0) rectangle (11,-1);
\fill[color=gray] (6,-1) rectangle (7,-3);
\draw (6, 0) grid (11,-1);
\draw (6,-1) grid (7,-3);

\fill[color=gray] (2,-5) rectangle (4,-6);
\fill[color=gray] (0,-6) rectangle (3,-7);
\fill[color=gray] (0,-7) rectangle (1,-9);
\draw (0, -5) grid (4,-6);
\draw (0,-6) grid (3,-7);
\draw (0,-7) grid (1,-9);

\fill[color=gray] (9,-5) rectangle (11,-6);
\fill[color=gray] (9,-6) rectangle (10,-8);
\fill[color=gray] (6,-7) rectangle (9,-8);
\draw (6,-5) grid (11,-6);
\draw (6,-6) grid (10,-8);

\fill[color=gray] (2,-10) rectangle (4,-11);
\fill[color=gray] (0,-11) rectangle (3,-12);
\fill[color=gray] (0,-12) rectangle (1,-14);
\draw (0,-10) grid (4,-11);
\draw (0,-11) grid (3,-12);
\draw (0,-12) grid (1,-14);

\fill[color=gray] (8,-10) rectangle (11,-11);
\fill[color=gray] (8,-11) rectangle (9,-13);
\fill[color=gray] (6,-12) rectangle (9,-13);
\draw (6,-10) grid (11,-11);
\draw (6,-11) grid (9,-13);

\fill[color=gray] (1,-15) rectangle (4,-16);
\fill[color=gray] (0,-16) rectangle (2,-17);
\fill[color=gray] (0,-17) rectangle (1,-19);
\draw (0,-15) grid (4,-16);
\draw (0,-16) grid (2,-17);
\draw (0,-17) grid (1,-19);

\fill[color=gray] (7,-15) rectangle (11,-16);
\fill[color=gray] (7,-16) rectangle (8,-18);
\fill[color=gray] (6,-17) rectangle (7,-18);
\draw (6, -15) grid (11,-16);
\draw (6,-16) grid (8,-18);

\end{tikzpicture}
\end{tabular}

\caption{Example of the partitions enumerated in
Theorem \ref{theorem_fu_tang_generalization} with $m=1$, $k=1$,
and $n=7$.}
\label{fu_tang_generalization_example}
\end{figure}

The rest of this paper is organized into two sections. Section 2 contains a
description of the technique used in the proofs, as well as the proof of Theorem
\ref{theorem_fu_tang_generalization}. The combinatorial proof of Theorem
\ref{main_result} is in section 3, as well as some commentary.

\section{Compositions and the partition perimeter}

A \textit{composition} of the integer $n \geq 1$ is a partition
$c = (c_1, c_2, \dots, c_\ell)$ of $n$ with the requirement that
$c_1 \geq c_2 \geq \cdots \geq c_\ell$ discarded. We borrow the same terminology and
tools for compositions from partitions when appropriate.

We may visually represent $c$ as a \textit{Young diagram}, which is a grid
of squares containing $c_i$ squares in the $i$th row from the top, aligned on the left.
The \textit{$m$-modular diagram} of $c$ for fixed $m\geq 2$ is a Young diagram where the
squares furthest to the right contain a value between 1 and $m$, and all others contain
$m$, having the property that the sum of the values in the $i$th row is $c_i$.
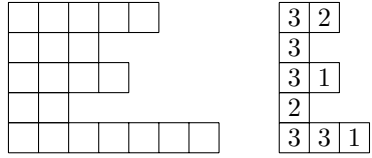
\begin{figure}[!ht]
\centering
\begin{tikzpicture}[scale=0.4]
\draw (0,0) grid (5,-1);
\draw (0,-1) grid (3,-2);
\draw (0,-2) grid (4,-3);
\draw (0,-3) grid (2,-4);
\draw (0,-4) grid (7,-5);

\draw (9,0) grid (11,-1);
\draw (9,-1) grid (10,-2);
\draw (9,-2) grid (11,-3);
\draw (9,-3) grid (10,-4);
\draw (9,-4) grid (12,-5);

\node at (9.5,-0.5) {3};
\node at (10.5,-0.5) {2};
\node at (9.5,-1.5) {3};
\node at (9.5,-2.5) {3};
\node at (10.5,-2.5) {1};
\node at (9.5,-3.5) {2};
\node at (9.5,-4.5) {3};
\node at (10.5,-4.5) {3};
\node at (11.5,-4.5) {1};
\end{tikzpicture}
\caption{The Young diagram and $3$-modular diagram of $c = (5, 3, 4, 2, 7)$.}
\end{figure}

For any partition $\lambda$, the \textit{conjugate} $\lambda'$ is a new partition
with $\lambda_i'$ defined to be the $i$th column of $\lambda$'s Young diagram, from the
left. Conjugation is an involution that
preserves the perimeter since $\lambda_1$ and $\ell(\lambda)$ are swapped.
\newline

\noindent \textbf{Remark.}
Conjugation is only defined for partitions, although an analogous operation
exists for compositions in general (See the work of Munagi such as \cite{Munagi} and
the references therein for recent work).

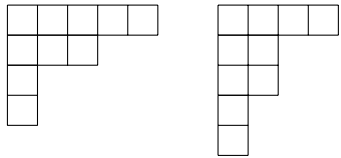
\begin{figure}[!ht]
\centering
\begin{tikzpicture}[scale=0.4]
\draw (0,0) grid (5,-1);
\draw (0,-1) grid (3, -2);
\draw (0,-2) grid (1, -4);
\draw (7,0) grid (11,-1);
\draw (7,-1) grid (9,-3);
\draw (7,-3) grid (8,-5);
\end{tikzpicture}
\caption{The conjugate partition pair $(5,3,1,1)$ and $(4,2,2,1,1)$.}
\end{figure}

The perimeter of a partition $\lambda$ can be equivalently defined as the number
of squares with an edge or corner on the bottom or right side of $\lambda$'s
Young diagram. Let
\begin{equation*}
c = (\lambda_1 - \lambda_2 + 1, \lambda_2 - \lambda_3 + 1, \dots,
\lambda_{\ell(\lambda)-1} - \lambda_{\ell(\lambda)} + 1, \lambda_{\ell(\lambda)}).
\end{equation*}

Clearly $c_i$ counts the number of squares from the $i$th row of $\lambda$'s Young
diagram that contribute to the perimeter, so $c$ is a composition with the properties
$\lvert c \rvert = \lambda_1 + \ell(\lambda) - 1$ and $\ell(c) = \ell(\lambda)$.
Indeed, this defines a bijection since each part of $\lambda$ can easily be recovered
using the formula
\begin{equation*}
\lambda_i = 1 + \sum_{j=i}^{\ell(c)} (c_j - 1).
\end{equation*}
Let $\mathcal{C}$ and $\mathcal{P}$ be the sets of all compositions and partitions,
respectively, and $\pi \colon \mathcal{C} \to \mathcal{P}$ be this bijection. We can
also interpret $\pi(c)$ as the partition formed by sliding the rows
of $c$ over so they overlap at exactly one square. 

\begin{figure}[!ht]
\centering
\begin{tikzpicture}[scale=0.4]
\fill[color=gray] (0,0) rectangle (1,-5);
\fill[color=gray] (1,0) rectangle (3,-1);
\fill[color=gray] (1,-2) rectangle (2,-3);
\fill[color=gray] (1,-3) rectangle (4,-4);
\fill[color=gray] (1,-4) rectangle (2,-5);
\draw (0,0) grid (1,-5);
\draw (1,0) grid (3,-1);
\draw (1,-2) grid (2,-3);
\draw (1,-3) grid (4,-4);
\draw (1,-4) grid (2,-5);
\fill[color=gray] (13,0) rectangle (16,-1);
\fill[color=gray] (13,-1) rectangle (14,-2);
\fill[color=gray] (12,-2) rectangle (14,-3);
\fill[color=gray] (9,-3) rectangle (13,-4);
\fill[color=gray] (8,-4) rectangle (10,-5);
\draw (8,0) grid (16, -1);
\draw (8,-1) grid (14, -3);
\draw (8,-3) grid (13,-4);
\draw (8,-4) grid (10,-5);
\draw[ultra thick, ->] (5,-2.5) -- (7, -2.5);
\node at (6,-3.75) {$p(c)$};
\end{tikzpicture}
\caption{An example of the bijection $\pi$ with $c = (3,1,2,4,2)$ and $\pi(c)=(8,6,6,5,2)$.}
\label{pi_figure}
\end{figure}
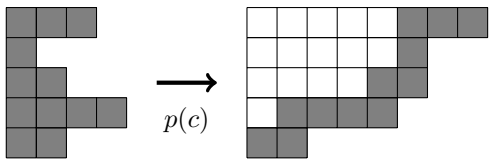

Using $\pi$ and these Young diagram interpretations,
we can encode information about the partition perimeter into a composition,
which allows us to then work entirely with compositions, only switching back with $\pi^{-1}$
when done.

For an example of how this technique may be used,
consider the following composition result of Munagi.

\begin{theorem}[Theorem 1.2 in \cite{Munagi2}]
\label{theorem_munagi}
For all $n, m \geq 1$ there are as many compositions of $n$ with parts $\equiv 1 \pmod{m}$
as there are compositions of $n + m - 1$ with parts $\geq m$.
\end{theorem}

Replacing $m$ with $m+1$, the composition parts from the first family become
under $\pi$ partition parts $\equiv 1 \pmod{m+1}$ since they all have the effect
of preserving the congruence of the part below. Deleting $m$ from the second family's
last part, which is guaranteed to be $\geq m+1$, brings the size down to $n$. 
Sliding the parts over as in Figure \ref{pi_figure}, we obtain
partitions with gaps between parts $\geq m$, and so we arrive at Theorem
\ref{theorem_fu_tang}.

The proof of Theorem \ref{theorem_fu_tang_generalization} is similar to this,
using the following refinement of Munagi's theorem from Huang.

\begin{theorem}[Theorem 1.6 in \cite{Huang}]
\label{theorem_huang}
For all $n, m \geq 1$, and $k \geq 0$, there are as many compositions of $n$
with exactly $k$ parts $\not\equiv 1 \pmod{m}$, each $> m$,
as there are compositions of $n+m-1$ with
exactly $k$ parts $< m$, each proceeded by a part  $\geq m$ and
followed by either a last part or a part $>m$.
\end{theorem}

\noindent\textit{Proof of Theorem \ref{theorem_fu_tang_generalization}}.
Proceed as above, replacing $m$ with $m+1$ now using Theorem
\ref{theorem_huang}.

For family 1), any composition parts $\not\equiv 1 \pmod{m+1}$ have the effect
under $\pi$ of changing from
the part below's congruence class, if there is one,
so each of the $k$ sets are formed from
a part $\not\equiv 1 \pmod{m+1}$, followed above by any amount of parts
$\equiv 1 \pmod{m+1}$.

For family 2), delete $m$ from the final part of each composition. Each of the $k$
parts $<m+1$ along with the guaranteed part above $\geq m+1$ and below $> m+1$
produce under $\pi$ a gap $\geq m$, $<m$, and then $>m$.

The omitted details are straight forward to verify.

\qed

\noindent \textbf{Remark.} Huang provides a formula (Theorem 1.7 in \cite{Huang})
that the partitions in Theorem \ref{theorem_fu_tang_generalization} inherit.
\newline

Using this technique in the other direction is possible too, yielding composition
results. Although this idea is not pursued in this work, we provide the following
simple example.

\begin{proposition}
\label{example_proposition}
For all $m \geq 2 $ and $n, k\geq 1$, there are as many compositions of $n$ with $k$
parts, all at most $m$ and last part
less than $m$, as there are compositions of $n$ with
$n - k + 1$ parts and no continuous sequence of $m-1$ or more $1$'s, excluding the
last part.
\end{proposition}

\noindent \textit{Proof.}
Partitions with perimeter $n$, $k$ parts, and gaps between parts $< m$
conjugate to give partitions
with perimeter $n$, $n-k+1$ parts since $k-1$ squares in the Young diagram
had a perimeter contributing square below them,
and parts repeating $< m$ times. Interpreting these families
through $\pi^{-1}$ gives the result.

\qed

\section{Combinatorial proof of Theorem \ref{main_result}}

Before proceeding, we need the following lemma.
Here we use the notation $\equiv R \;(\textrm{mod}\; m)$ to mean
congruent to some $r \in R$ modulo $m$.

\begin{lemma} Let $m \geq 2$, $n \geq 1$, and
$R \subseteq \{1, 2, \dots, m-1\}$
be non-empty. There are as many compositions of $n$ with parts
$\equiv R \;(\textrm{mod}\; m)$ as there are compositions
of $n$ with parts in $R \cup \{m\}$, with last part in $R$.
\label{congruence_lemma}
\end{lemma}

%We offer a proof of Lemma \ref{congruence_lemma} in two flavors, analytic and bijective.
%The latter will be needed for the proof of Theorem \ref{main_result}.
%\newline
%
%\textit{Analytic proof of Lemma \ref{congruence_lemma}.}
%Write $R = \{r_1, \dots, r_k\}$. By summing over compositions
%with number of parts $i$, the generating function for compositions
%with parts $\equiv R \;(\textrm{mod}\; m)$ is
%\begin{align*}&\sum_{i=1}^\infty \left( \frac{q^{r_1} + \cdots + q^{r_k}}{1-q^m} \right)^i \\
%=\quad & \frac{1}{1 - (q^{r_1} + \cdots + q^{r_k})/(1-q^m)} - 1 \\
%=\quad & \frac{q^{r_1} + \cdots + q^{r_k}}{1- (q^{r_1} + \cdots + q^{r_k} + q^m)} \\
%=\quad & (q^{r_1} + \cdots + q^{r_k})\sum_{i=0}^\infty
%(q^{r_1} + \cdots + q^{r_k} + q^m)^i
%\end{align*}
%which after some standard manipulations gives us the generating
%function for compositions with parts in $R \cup \{m\}$ and last
%part in $R$, again summed by number of parts $i$.
%
%\qed

\textit{Proof.}
Let $c$ be a composition with parts $\equiv R \;(\textrm{mod}\; m)$.
Writing $c$ as an $m$-modular diagram, rotate each part one quarter
turn clockwise, leaving the $m$-modular diagram of a composition
with parts in $R\cup \{m\}$ and last part in $R$. Since the numbers
were unchanged, this resulting composition has the same size.

For the inverse, given a composition $c$ with parts
in $R\cup\{m\}$ and last part in $R$, iterate upwards in the $m$-modular
diagram starting from
$c$'s last part. Collect the square not containing $m$ and then all
squares containing $m$ directly above, then rotate counter-clockwise
to form a part $\equiv R \;(\mathrm{mod}\; m)$. Repeat until
all squares are exhausted.

\qed

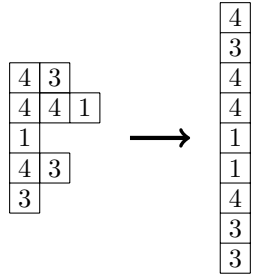
\begin{figure}[!ht]
\centering
\begin{tikzpicture}[scale=0.4]
\draw[step=1] (0,0) grid (1,-5);
\draw[step=1] (1,0) grid (2,-1);
\draw[step=1] (1,-1) grid (3,-2);
\draw[step=1] (1,-3) grid (2,-4);
\node at (0.5, -0.5) {4};
\node at (0.5, -1.5) {4};
\node at (0.5, -2.5) {1};
\node at (0.5, -3.5) {4};
\node at (0.5, -4.5) {3};
\node at (1.5, -0.5) {3};
\node at (1.5, -1.5) {4};
\node at (1.5, -3.5) {3};
\node at (2.5, -1.5) {1};
\draw[ultra thick, ->] (4, -2.5) -- (6, -2.5);
\draw[step=1] (7,2) grid (8,-7);
\node at (7.5, 1.5) {4};
\node at (7.5, 0.5) {3};
\node at (7.5, -0.5) {4};
\node at (7.5, -1.5) {4};
\node at (7.5, -2.5) {1};
\node at (7.5, -3.5) {1};
\node at (7.5, -4.5) {4};
\node at (7.5, -5.5) {3};
\node at (7.5, -6.5) {3};
\end{tikzpicture}
\caption{ With $m=4$ and $R=\{1,3\}$,
the composition $(7,9,1,7,3)$ under this bijection maps to
$(4,3,4,4,1,1,4,3,3)$.}
\end{figure}

\textit{Combinatorial proof of Theorem \ref{main_result}.}
Fix the integers $m$ and $n$. The proof is built using a series of simple bijections
and one injection between the following six families of compositions and partitions.

\begin{itemize}
\item[\textbf{1)}] Partitions with perimeter $n$ and parts $\not\equiv 0 \pmod{m}$ become
\item[\textbf{2)}] compositions $c$ of $n$ with the property that
\begin{equation}\label{equation_m_regular}\tag{$\ast$}
1 + \sum_{j=i}^{\ell(c)} (c_j - 1) \not\equiv 0\;
  (\mathrm{mod}\; m)\; \mathrm{for\; all}
  \; 1 \leq i \leq \ell(c)
\end{equation}
through $\pi^{-1}$. The injection $\varphi$ defined below sends these to
\item[\textbf{3)}] compositions of $n$ with parts $\not \equiv 0 \pmod{m}$. Next, using
Lemma \ref{congruence_lemma} with $R=\{1,\dots,m-1\}$, these become
\item[\textbf{4)}] compositions of $n$ with all parts $\leq m$ and last part $< m$. Applying
$\pi$, we then get
\item[\textbf{5)}] partitions with perimeter $n$ and gaps between parts $< m$ and then
finally
\item[\textbf{6)}] partitions with perimeter $n$ and parts repeating $<m$ times, by
conjugation.
\end{itemize}

The injection $\varphi$ is defined as follows.
Let $c$ be a composition
with the property (\ref{equation_m_regular}). From $c$ we
construct $\varphi(c) = d$ using the following iterative algorithm,
starting at $c$'s last part. The algorithm is designed to preserve parts in $c$ that are
already $\not\equiv 0\; (\mathrm{mod}\; m)$, and transform the others in such a
way that a part in $\pi(d)$, up to a small technical detail,
becomes $\equiv 0 \pmod{m}$, violating
(\ref{equation_m_regular}) with each occurrence.
Here, we use the notation $()$ to mean the \textit{empty composition}, which has no
parts and size 0.
\newline

\begin{algorithmic}[1]
\STATE Initialize with $i = \ell(c)$,  $j = 1$, and $d = ()$
\WHILE{$i \neq 0$}
\STATE Let $j = j + c_i - 1$.
\IF{$c_i \not\equiv 0\; (\mathrm{mod}\; m)$}
\STATE Append $c_i$ to the start of $d$ ($c_i$ becomes the new first part of $d$)
\ELSE
\STATE Let $r$ be the remainder after dividing $j$ by $m$
\STATE Append the part $m-r$ to the start of $d$
\STATE Append the part $c_i - (m-r)$ to the start of $d$
\ENDIF
\STATE Let $i = i - 1$
\ENDWHILE
\end{algorithmic}
\; \newline
From this construction, that $\lvert c \rvert = \lvert d \rvert$ is immediate. Moreover,
from (\ref{equation_m_regular}), $j \not\equiv 0\; (\mathrm{mod}\; m)$ will hold at
every step, so in the case of lines 7 through 9, $r$ will never be 0, which implies
that the parts $m-r$ and $c_i - (m-r)$ are $\not\equiv 0\; (\mathrm{mod}\; m)$ as well.
Therefore each part of $d$ is $\not\equiv 0\; (\mathrm{mod}\; m)$.

\begin{figure}[!ht]
\centering
\begin{tikzpicture}[scale=0.4]
\fill[color=gray] (0,0) rectangle (6,-1);
\fill[color=gray] (0,-3) rectangle (3,-4);
\draw (0,0) grid (6,-1);
\draw (0,-1) grid (2,-2);
\draw (0,-2) grid (4,-3);
\draw (0,-3) grid (3,-4);
\draw (0,-4) grid (2,-5);
\draw[ultra thick, ->] (7,-2.5) -- (9,-2.5);
\node at (8,-3.75) {$\varphi(c)$};
\fill[color=gray] (10,1) rectangle (14,0);
\fill[color=gray] (10,0) rectangle (12,-1);
\fill[color=gray] (10,-3) rectangle (11,-4);
\fill[color=gray] (10,-4) rectangle (12,-5);
\draw (10,1) grid (14,0);
\draw (10,0) grid (12,-1);
\draw (10,-1) grid (12,-2);
\draw (10,-2) grid (14,-3);
\draw (10,-3) grid (11,-4);
\draw (10,-4) grid (12,-5);
\draw (10,-5) grid (12,-6);
\end{tikzpicture}
\caption{With $m=3$, $\varphi$ sends $c = (6,2,4,3,2)$ to $d=(4,2,2,4,1,2,2)$.}
\end{figure}
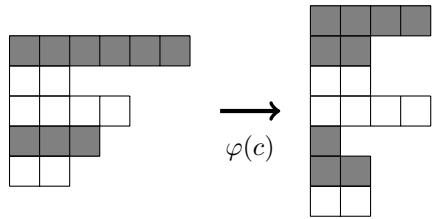

Suppose that some composition $d$ is in the image
of $\varphi$. Like in the
algorithm, we iterate $i$ from $d$'s last part to the top, with the goal
of uniquely determining the parts of the preimage. Start with $c=()$ and $j=1$,
and increment $j$ by $d_i - 1$ at each step.

If $j \not\equiv 0 \;(\mathrm{mod}\; m)$
then we must be in the case of line $5$, so append $d_i$ to the top of $c$.
Otherwise, we are in the case of lines 7 through 9 and so we must append
$d_i + d_{i-1}$ to the start of $c$. We must then subtract 1 from $i$ an additional time,
and also add 1 to $j$ since the action of splitting the composition part into two in $c$
has the effect of decreasing $j$ by 1 in $d$.
Continue until all parts of $d$ are exhausted, which leaves a uniquely determined $c$.
This shows injectivity which completes the proof.

\qed

\subsection{Closing Remarks}
If $n > m > 2$ in Theorem \ref{main_result} the inequality
becomes strict since the compositions of size $n$
\begin{equation*}
(m-1, 2, \hspace{-0.45cm}
\underbrace{1, \dots, 1}_{n-m-1 \;\mathrm{many\; 1's}}
\hspace{-0.45cm})
\end{equation*}
cannot be in the image of $\varphi$. This fact can also be obtained
from the original proof. If $n \leq m$ both
families of partitions coincide giving an equality. Also, if $m=2$, $\varphi$ becomes
the identity map, which provides another proof of Straub's theorem.

Many papers that treat compositions have noted that
the following families of compositions are all enumerated by the $n$th Fibonacci
number $F_n$. Recall that $F_0 = 0, F_1 = 1$, and $F_n = F_{n-1} + F_{n-2}$ for
any $n \geq 2$.

\begin{itemize}
\item[1)] Compositions of $n$ with odd parts.
\item[2)] Compositions of $n$ with parts 1 and 2, and last part 1.
\item[3)] Compositions of $n+1$ with parts greater than 1.
\end{itemize}

Lemma \ref{congruence_lemma} with $m=2$ and $R = \{1\}$ gives a bijection from
1) to 2), and by identifying 3) with compositions of $n$ where only the last
part may be 1, Proposition \ref{example_proposition} with $m=2$ gives a bijection from 2)
to 3).

\bibliography{references.bib}
\bibliographystyle{amsplain}

\end{document}